\documentclass[12pt]{article}

\usepackage{amsmath,amsthm,amsfonts,amssymb,amscd}
\def\qed{{\unskip\nobreak\hfil\penalty50
\hskip2em\hbox{}\nobreak\hfil$\square$
\parfillskip=0pt \finalhyphendemerits=0\par}\medskip}

\def\Ad{{\mathrm {Ad}}}

\def\Tr{{\mathrm {Tr}}}

\def\dim{{\mathrm {dim}}}

\def\Car{{\mathrm{CAR}}}
\def\Ga{\Gamma}

\def\Mob{{\rm\textsf{M\"ob}}}

\def\Ad{{\mathrm {Ad}}}

\def\Tr{{\mathrm {Tr}}}

\def\dim{{\mathrm {dim}}}

\def\Ga{\Gamma}

\newtheorem{theorem}{Theorem}[section]
\newtheorem{lemma}[theorem]{Lemma}
\newtheorem{conjecture}[theorem]{Conjecture}
\newtheorem{corollary}[theorem]{Corollary}
\newtheorem{definition}[theorem]{Definition}

\newtheorem{proposition}[theorem]{Proposition}
\newtheorem{remark}[theorem]{Remark}

\def\Mob{{\rm\textsf{M\"ob}}}

\def\res{\!\restriction\!}
\def\A{{\cal A}}
\def\g{{\mathbf g}}
\def\B{{\cal B}}
\def\C{{\cal C}}

\def\I{{\cal I}}

\def\PI{{\cal PI}}

\def\H{{\cal H}}
\renewcommand{\qed}{\ \hfill $\blacksquare$}

\newcommand{\bdef}{\begin{definition}}
\newcommand{\blem}{\begin{lemma}}
\newcommand{\bprop}{\begin{proposition}}
\newcommand{\bthm}{\begin{theorem}}
\newcommand{\bcor}{\begin{corollary}}
\newcommand{\bconj}{\begin{conjecture}}
\newcommand{\ben}{\begin{equation}}
\newcommand{\een}{\end{equation}}

\newcommand{\ede}{\end{definition}}
\newcommand{\elem}{\end{lemma}}
\newcommand{\eprop}{\end{proposition}}
\newcommand{\ethm}{\end{theorem}}
\newcommand{\ecor}{\end{corollary}}
\newcommand{\econj}{\end{conjecture}}
\newcommand{\brem}{\begin{remark}}
\newcommand{\erem}{\end{remark}}

\newcommand{\ba}{\begin{array}}
\newcommand{\ea}{\end{array}}
\newcommand{\bea}{\begin{eqnarray}}

\title{\huge Relative Entropy in CFT}
\author{
{\sc Roberto Longo}\footnote{Supported by ERC Advanced Grant 669240 QUEST ``Quantum Algebraic Structures and Models'', MIUR FARE R16X5RB55W  QUEST-NET, GNAMPA-INdAM and Alexander von Humboldt Foundation.}
\\
Department of Mathematics\\
University of Rome Tor Vergata\\
Via della Ricerca Scientifica, 1, 00133 Roma, Italy\\
E-mail: {\tt longo@mat.uniroma2.it}\\
{}
\\
{\sc Feng Xu}$^{**}$\\
Department of Mathematics\\
University of California at Riverside\\
Riverside, CA 92521\\
E-mail: {\tt xufeng@math.ucr.edu}}

\begin{document}
\date{}
\maketitle

\begin{abstract}
By using Araki's relative entropy, Lieb's convexity and the theory of
singular integrals, we compute the mutual information associated with free
fermions, and we deduce many results about entropies for chiral CFT's
which are embedded into free fermions, and their extensions. Such relative entropies in CFT are
here computed explicitly for the first time in a mathematical rigorous way.  Our results
agree with previous computations by physicists based on heuristic arguments;
in addition we uncover a surprising connection with the theory of
subfactors, in particular by showing that a certain duality, which is
argued to be true on physical grounds, is in fact violated if the
global dimension of the conformal net is greater than $1.$

\end{abstract}

\newpage

\section{Introduction}
In the last few years there has been an enormous amount of work by physicists
concerning
entanglement entropies in QFT, motivated by the connections
with condensed matter physics, black holes, etc.; see the references in
\cite{Hollands} for a partial list of references.  However, some very basic
mathematical questions remain open. For example, most of the
entropies computed in the physics literature are infinite, so the
singularity structures, and sometimes the cut off independent
quantities,  are of most interest. Often, the mutual information
is argued to be finite based on heuristic physical arguments, and
one can derive the singularities of the entropies from the mutual information by
taking singular limits. But it is not even clear that such mutual
information, which is well defined as a special case of
Araki's relative entropy, is indeed finite.

In this paper we begin to
address some of these fundamental mathematical questions motivated
by the physicists' work on entropy.  For related works, see
\cite{Hollands} and \cite{OT}.  Unlike the main focus in
\cite{Hollands}, the mutual information considered in our paper can
be computed explicitly in many cases and satisfies many conditions,
but not all, proposed by physicists such as those  in
\cite{Casc}. Our project is strongly motivated by Edward Witten's
questions, in particular his question to  make physicists' entropy
computations rigorous. In this paper we focus on the Chiral CFT in
two dimensions, where the results we have obtained are most explicit
and have interesting connections to subfactor theory, even though
some of our results, such as Theorem \ref{main2}, do not depend on
conformal symmetries and apply to more general QFT.  The main
results are:

1) Theorem \ref{main1}: Exact computation of the mutual information (through the relative entropy as defined by Araki for general states on von Neumann algebras)  for free fermions.  Note that this was not even known to be finite, for example the main quantity defined in \cite{Hollands} is smaller and does not seem to
 verify the conditions in the physical literature.
Our proof uses Lieb's convexity and the theory of singular integrals; to the best of our knowledge, by Theorem \ref{main1},  Theorem \ref{soft2} and examples in
Section \ref{examples}, this is the first time
 that such relative entropies are computed in a mathematical rigorous way. The results verify earlier computations by physicists based on heuristic arguments, such as P. Calabrese and J. Cardy  in \cite{CC} and H. Casini and  M. Huerta in \cite{Casfree}.

In particular, for the free chiral net $\A_r$ associated with $r$ fermions, and two intervals
$A=(a_1, b_1)$, $B=(a_2, b_2)$ of the real line, where $b_1< a_2$, the mutual information
associated with $A,B$ is
$$F (A,B)= -\frac{r}{6}\ln \eta \ ,$$
where $\eta = \frac{(b_1-a_2)(b_2-a_1)}{(b_1-a_1)(b_2-a_2)}$ is the
cross ratio of $A,B$, $0< \eta <1$.

2) It follows from 1) and the monotonicity of the relative entropy that
any chiral CFT in two dimensions that embeds into free fermions, and
their finite index extensions,
 verify most of the conditions (not all, see Section \ref{fail})
discussed for example in \cite{Casc}, see Theorem \ref{soft1}. This
includes a large family of chiral CFTs. Much more can be obtained if
the embedding has finite index as in Theorem \ref{soft2}. In this
case, we also verify a proposal (cf. (1) of Theorem \ref{main2} ) in
\cite{Casc} about an entropy formula related to  a derivation of the
$c$ theorem. Theorem \ref{soft2} also connects relative entropy and
index of subfactors in an interesting and unexpected way. There is
one bit of surprise: it is usually  postulated that the mutual
information of a pure state such as vacuum state for complementary
regions should be the same. But in the Chiral case this is not true,
and the violation is measured by global dimension of the chiral CFT
as will be seen in Section \ref{fail}.

The physical meaning of the last part of (2) is not clear to us. The
violation, which is in some sense proportional to the logarithm of
global index, also turns out to be what is  called topological
entanglement entropy (cf. Remark \ref{kitaevr}). In \cite{IW} the
authors discuss chiral theories where  entanglement entropy cannot
be defined with the expected properties due to anomalies. The
relation to our work is not clear. On the other hand, when
considering full CFT, one does have global dimension equal to $1$,
and it remains an interesting question to investigate entropies in
the full CFT framework.

The rest of this paper is structured as follows. After a preliminary
section on von Neumann entropy, Araki's relative entropy, graded
nets and subnets, we consider the computation of mutual information
in \S\ref{fermion}. In \S\ref{descent}, we derive many of the
properties of the mutual information in the vacuum state for all Chiral CFT
which are embedded into free fermions, and their extensions, from the
results of \S\ref{fermion}. In the last section we supply with two
families of chiral CFT where our main results apply.

 \section{Preliminaries}\label{prelim}

\subsection{Entropy and relative entropy}

von Neumann entropy is the quantity associated with a density matrix $\rho$ on a
Hilbert space $\H$ by
\[
S(\rho) = -\Tr(\rho \log \rho)\ .
\]
von Neumann entropy can be viewed as a measure of the lack of information about a
system to which one has ascribed the state. This interpretation is
in accord for instance with the facts that $S(\rho)\ge0$ and that a
pure state $\rho = |\Psi \rangle \langle \Psi|$ has vanishing von
Neumann entropy.

A related notion is that of the relative entropy. It is defined for
two density matrices $\rho, \rho'$ by \ben  S(\rho,\rho') = \Tr(\rho
\log \rho - \rho \log \rho') \ . \een  Like $S(\rho)$,
$S(\rho,\rho')$ is non-negative, and can be infinite.

A generalization of the relative entropy in the context of von Neumann algebras of
arbitrary type was found by Araki~\cite{arakir} and is formulated
using modular theory. Given two faithful, normal states $\omega,
\omega'$ on a von Neumann algebra $\A$ in standard form, we choose
the vector representatives in the natural cone $\mathcal{P}^\sharp$,
called $|\Omega \rangle,|\Omega' \rangle$ . The anti-linear opeartor
 $S_{\omega,\omega'} a|\Omega' \rangle = a^*
|\Omega \rangle$, $a\in \A$, is closable
and one considers again the polar decomposition of its closure
$\bar S_{\omega, \omega'} = J \Delta_{\omega,\omega'}^{1/2}$ . Here $J$ is the
modular conjugation of $\A$ associated with $\mathcal{P}^\sharp$ and $\Delta_{\omega,\omega'} =
 S^*_{\omega, \omega'}\bar S_{\omega, \omega'}$ is the relative modular operator w.r.t. $|\Omega \rangle,|\Omega' \rangle$.
Of course, if $\omega = \omega'$ then $\Delta_{\omega} = \Delta_{\omega,\omega'}$ is the usual modular operator.

A related
object is the Connes cocycle (Radon-Nikodym derivative) defined as
$[D\omega : D\omega']_t = \Delta_{\omega, \psi}^{it} \Delta_{\psi,
\omega'}^{it} \in \A$, where $\psi$ is an arbitrary auxiliary
faithful normal state on $\A'$ . \bdef
 The relative entropy w.r.t. $\omega$ and $\omega'$ is defined by
  \ben\label{drel1}
S(\omega, \omega') = \langle \Omega | \log \Delta_{\omega, \omega'}
\ \Omega\rangle  = \lim_{t \to 0} \frac{\omega([D\omega:D\omega']_t
- 1)}{it} \ , \een $S$ is extended to positive linear functionals that are
not necessarily normalized by the formula
$S(\lambda\omega,\lambda'\omega')=\lambda S(\omega,\omega') +
\lambda \log (\lambda/\lambda')$, where $\lambda, \lambda'>0$ and
$\omega,\omega'$ are normalized. If $\omega'$ is not normal, then
one sets $S(\omega, \omega') = \infty$. \ede

 For a type I algebra $\A = \B(\H)$, states $\omega, \omega'$
correspond to density matrices $\rho, \rho'$. The square root of the relative modular
operator $\Delta_{\omega,\omega'}^{1/2}$ corresponds to $\rho^{1/2}
\otimes \rho^{\prime -1/2}$ in the standard representation of $\A$ on $\H
\otimes \bar \H$; namely $\H
\otimes \bar \H$ is identified with the Hiilbert-Schmidt operators $HS(\H)$ with the left/right multiplication
of $\A$/$\A'$.
 In this representation, $\omega$ corresponds to
the vector state $|\Omega \rangle =\rho^{1/2} \in \H \otimes
\bar{\H}$, and the abstract definition of the relative entropy
in~\eqref{drel1} becomes \ben \langle \Omega | \log \Delta_{\omega,
\omega'} \, \Omega \rangle = \Tr_{\H} \rho^{\frac12}\left(\log\rho
\otimes 1 - 1 \otimes \log \rho'\right)\rho^{\frac12} = \Tr_\H(\rho
\log \rho - \rho \log \rho') \ . \een

\medskip

As another example, let us consider a bi-partite system with Hilbert
space $\H_A \otimes \H_B$ and observable algebra $\A = \B(\H_A)
\otimes \B(\H_B)$. A  normal state $\omega_{AB}$ on $\A$ corresponds to a
density matrix $\rho_{A B}$.  One calls $\rho_A = \Tr_{\H_B}
\rho_{AB}$ the ``reduced density matrix'', which defines a state
$\omega_A$ on $\B(\H_A)$ (and similarly for system $B$).
 The mutual
information is given in our example system by
\ben\label{remu}
S(\rho_{AB},\rho_A\otimes \rho_B ) = S(\rho_A) +S(\rho_B) -
S(\rho_{A B}) \ . \een

For tri-partite system with Hilbert space $\H_A \otimes \H_B \otimes
\H_C$ and observable algebra $\A = \B(\H_A) \otimes \B(\H_B)\otimes
\B(\H_C)$, we have the following strong subadditivity (cf.
\cite{Lieb}):
\ben\label{strong} S(\rho_{AB}) + S(\rho_{A
C})-S(\rho_A)-S(\rho_{A B C})\geq 0\ . \een
A list of properties of relative entropies that will be used later
can be found in \cite{OP} (cf. Th. 5.3, Th. 5.15 and Cor. 5.12
\cite{OP}):

\bthm\label{515} (1) Let $M$ be a von Neumann algebra and $M_1$ a
von Neumann subalgebra of M. Assume that there exists a faithful
normal conditional expectation $E$ of $M $onto $M_1$. If $\psi$ and
$\omega$  are states of $M_1$ and $M$, respectively, then $S(\omega,
\psi\cdot E) = S(\omega\res M_1, \psi) + S(\omega, \omega\cdot E);$
\par

(2) Let be $M_i$ an increasing net of von Neumann subalgebras of $
M$ with the property $ (\bigcup_i M_i)''=M$. Then $S(\omega_1\res
M_i, \omega_2\res M_i)$ converges to $ S(\omega_1,\omega_2)$ where
$\omega_1, \omega_2$ are two normal states on $M$; \par

(3) Let $\omega$ and  $\omega_1$ be two normal states on a von
Neumann algebra $M$. If $\omega_1\geq \mu\omega,  $ then $S(\omega,
\omega_1) \leq \ln \mu^{-1}$;

(4) Let $\omega$ and  $\phi$ be two normal  states on a von Neumann
algebra $M$, and denote by  $\omega_1$ and  $\phi_1$ the
restrictions of   $\omega$ and  $\phi$ to  a von Neumann subalgebra
$M_1\subset M$ respectively. Then $S(\omega_1, \phi_1)\leq S(\omega,
\phi)$.
\ethm
For type  $\mathrm{III}$  factors, the von Neumann entropy is always
infinite, but we shall see that in many cases mutual information is
finite. By taking singular limits, we can also explore the
singularities of von Neumann entropy from mutual information (cf.
\ref{soft2} for an example) which is important from physicists'
point of view. The formal properties of von Neumann entropies are
useful in proving properties of mutual information, see the proof of
Th. \ref{soft1}.

\subsection{Graded  nets and subnets}

This section is contained in \cite{Carp}. We refer to \cite{Carp}
for more details and proofs.

 We shall denote by $\Mob$ the M\"obius group, which is
isomorphic to $SL(2,\mathbb R)/\mathbb Z_2$ and acts naturally and
faithfully on the circle $S^1$.

By an interval of $S^1$ we mean, as usual, a non-empty, non-dense,
open, connected subset of $S^1$ and we denote by $\I$ the set of all
intervals. If $I\in\I$, then also $I'\in\I$ where $I'$ is the
interior of the complement of $I$.  Intervals are disjoint if their
closure are disjoint. We will denote by $\PI$ the set which consists
of disjoint union of intervals.

A {\it  net $\A$ of von Neumann algebras on $S^1$} is a map
\[
I\in\I\mapsto\A(I)
\]
from the set of intervals to the set of von Neumann algebras on a
(fixed) Hilbert space $\H$ which verifies the \emph{isotony
property}:
\[
I_1\subset I_2\Rightarrow  \A(I_1)\subset\A(I_2)
\]
where $I_1 , I_2\in\I$.

A \emph{M\"obius covariant net} $\A$ of von Neumann algebras on
$S^1$ is a net of von Neumann algebras on $S^1$ such that the
following properties $1-4$ hold:
\begin{description}
\item[\textnormal{\textsc{1. M\"obius covariance}}:] {\it There is a
strongly continuous unitary representation $U$ of} $\Mob$ {\it on
$\H$ such that}
\[
U(g)\A(I)U(g)^*=\A({g}I)\ , \qquad g\in \Mob,\ I\in\I \ .
\]
\end{description}
\begin{description}
\item[$\textnormal{\textsc{2. Positivity of the energy}}:$]
{\it The generator of the rotation one-para\-meter subgroup
$\theta\mapsto U({\rm rot}(\theta))$ (conformal Hamiltonian) is
positive}, na\-me\-ly $U$ is a positive energy representation.

\end{description}
\begin{description}
\item[\textnormal{\textsc{3. Existence and uniqueness of the vacuum}}:]
{\it There exists a unit $U$-invariant vector $\Omega$ (vacuum
vector), unique up to a phase, and $\Omega$ is cyclic for the
von~Neumann algebra $\vee_{I\in\I}\A(I)$}.
\end{description}

A $\mathbb Z_2$-grading on $\A$ is an involutive  automorphism $\g
=\Ad \Gamma$ of $\A$, such that $\Gamma^2=1,$  $\Gamma \Omega =
\Omega, \Gamma\A(I) \Gamma = \A(I)$ for all $ I.$

Given the grading $\g$, an element $x$ of $\A$ such that $\g(x)=\pm
x$ is called homogeneous, indeed a Bose or Fermi element according
to the $\pm$ alternative, or simply even or odd elements. We shall
say that the degree $\partial x$ of the homogeneous element $x$ is
$0$ in the Bose case and $1$ in the Fermi case.

A  {\em M\"obius covariant graded net} $\A$ on $S^1$ is a $\mathbb
Z_2$-graded M\"obius covariant net satisfying graded locality,
namely a  M\"obius covariant net of von Neumann algebras on $S^1$
such that the following holds:
\begin{description}
\item[\textnormal{\textsc{4. Graded locality}}:]
{\it There exists a grading automorphism $\g$ of $\A$ such that,  if
$I_1$ and $I_2$ are disjoint intervals,}
\[
[x,y] = 0,\quad x\in\A(I_1), y\in\A(I_2) \ .
\]
\end{description}
Here $[x,y]$ is the graded commutator with respect to the grading
automorphism $\g$ defined as follows: if $x,y$ are homogeneous then
\[
[x,y]\equiv xy - (-1)^{\partial x \cdot \partial y}yx
\]
and, for the general elements $x,y$, it is  extended by linearity. When
the grading is trivial, i.e., when $\Gamma=1,$ we shall refer to
$\A$ as a {\it local net}.

Note the \emph{Bose subnet} $\A_b$, namely the $\g$-fixed point
subnet $\A^\g$ of degree zero elements, is local.

Moreover, setting
\[
Z\equiv \frac{1 - i\Ga}{1 - i} \ ,
\]
we have that the unitary $Z$ fixes $\Omega$ and
\[
\A(I')\subset Z\A(I)'Z^*
\]
(twisted locality w.r.t. $Z$).
\medskip

\begin{theorem}\label{Reeh-Schlieder} Let $\A$ be a M\"obius covariant Fermi net on
$S^1$.  Then $\Omega$ is cyclic and separating for each von Neumann
algebra $\A(I)$, $I\in\I$.
 \end{theorem}

If $I\in \I$, we shall denote by $\Lambda_I$ the one parameter
subgroup of $\Mob$ of ``dilation associated with $I$\ \!''.

\begin{theorem}
Let $I\in\I$ and $\Delta_I$, $J_I$ be the modular operator and the
modular conjugation of $(\A(I),\Omega)$. Then we have:

$(i)$:
\begin{equation}\label{BW}
\Delta_{I}^{it} = U(\Lambda_I(-2\pi t)), \ t\in\mathbb R,
\end{equation}

$(ii)$: $U$ extends to an (anti-)unitary representation of {\rm
$\Mob\ltimes\mathbb Z_2$} determined by
\[
U(r_I)=ZJ_I,\ I\in\I,
\]
acting covariantly on $\A$, namely
\[
U(g)\A(I)U(g)^*=\A(\dot{g}I)\quad g\in\text{\rm
$\Mob$}\ltimes\mathbb Z_2\ I\in \I\ .
\]
Here $r_I:S^1\to S^1$ is the reflection mapping $I$ onto $I'$.
\end{theorem}
\begin{corollary} {\em (Additivity)}
  Let $I$ and $I_i$ be intervals with $I\subset\cup_i I_i$. Then
  $\A(I)\subset\vee_i\A(I_i)$.
\end{corollary}
\medskip

\begin{theorem}
For every  $I\in\I$, we have:
\[
\A(I')=Z\A(I)'Z^* \ .
\]
\end{theorem}
In the following corollary, the grading and the graded commutator is
considered on $B(\H)$ w.r.t. $\Ad\Ga$.
\begin{corollary}\label{carp2}
$\A(I')= \big\{x\in B(\H):\ [x,y]=0\ \forall y\in\A(I)\big\}$.
\end{corollary}

Let now $G$ be a simply connected  compact Lie group. By Th. 3.2 of
\cite{FG}, the vacuum positive energy representation of the loop
group $LG$ (cf. \cite{PS}) at level $k$ gives rise to an irreducible
local net denoted by {\it ${\A}_{G_k}$}. By Th. 3.3 of \cite{FG},
every irreducible positive energy representation of the loop group
$LG$ at level $k$ gives rise to  an irreducible covariant
representation of ${\A}_{G_k}$. When no confusion arises we will
write ${\A}_{G_k}$ simply as $G_k$ as in the last section
\ref{examples}.

\par

Next we  recall some definitions from \cite{KLM} . Recall that
${\I}$ denotes the set of intervals of $S^1$. Let $I_1, I_2\in
{\I}$. We say that $I_1, I_2$ are disjoint if $\bar I_1\cap \bar
I_2=\emptyset$, where $\bar I$ is the closure of $I$ in $S^1$.
 Denote by ${\I}_2$ the set of unions of
disjoint 2 elements in ${\I}$. Let ${\A}$ be a graded M\"{o}bius
covariant net. For $E=I_1\cup I_2\in{\I}_2$, let $I_3\cup I_4$ be
the interior of the complement of $I_1\cup I_2$ in $S^1$ where $I_3,
I_4$ are disjoint intervals. Let
$$
{\A}(E):= \A(I_1)\vee \A(I_2), \quad \hat {\A}(E):= (\A(I_3)\vee
\A(I_4))'.
$$ Note that ${\A}(E) \subset  Z \hat{\A}(E)Z^{-1}$, and its index
will be denoted by  $\mu_{{\A}}$ and is called the $\mu$-index of
${\A}$ or global index of $\A$. This generalizes the usual
$\mu$-index of ${\A}$ when $\A$ is local.

Let $\A$ be a  graded M\"{o}bius net. By a {\it  M\"{o}bius subnet}
(cf. \cite{L1}) we shall mean a map
\[
I\in\I\to\B(I)\subset \A(I)
\]
that associates to each interval $I\in \I$ a von Neumann subalgebra
$\B(I)$ of $\A(I)$, which is isotonic
\[
\B(I_1)\subset \A(I_2), I_1\subset I_2,
\]
and   M\"{o}bius covariant with respect to the representation $U$,
namely
\[
U(g) \B(I) U(g)^*= \B(gI)
\] for all $g\in \Mob$ and $I\in \I$, and we also require  that
 $\Ad \Gamma$ preserves $\B$ as a set.  Note that by Lemma 13 of \cite{L1} for each $I\in
\I$ there exists a conditional expectation $E_I: \A(I)\rightarrow
\B(I)$ such that $E_I$ preserves the vector state given by the
vacuum of $\A$. Let $P$ be the projection onto the closed subspace
spanned by $\B(I)\Omega.$
\begin{definition}\label{indexd}
Let $\A$ be a  graded M\"{o}bius covariant net and  $\B\subset \A $
a subnet. We say   $\B\subset \A $  is of finite index if
$\B(I)\subset \A(I)$is of finite index  for some (and hence all)
interval $I$. The index will be denoted by $[\A:\B].$
\end{definition}
Assume that  $\B\subset \A$  has finite index and
$[\A:\B]=\lambda^{-1}.$ Let $I_1$ and $I_2$ be two intervals
obtained from an interval $I$ by removing an interior point, and let
$J_1\subset I_2.$ By \cite{LR} there are isometries $w_1\in \A(I_1),
v_1\in \B(I_1)$ such that $a=\lambda^{-1}E_I (aw_1^*) w_1, \forall
a\in \A(I).$ Let $e_1=w_1w_1^*$. Then
$$
Pe_1P=\lambda P, \quad e_1 v_1v_1^* e_1=\lambda e_1,\quad  \lambda^{-1} v_1^*
e_1 v_1= 1 \ .
$$
Similarly we have  $w_2\in \A(J_1), v_2\in \B(J_1)$ and
$e_2=w_2w_2^*$ which verify same relations as above. $e_1\in
\A(I_1), e_2\in \A(J_1)$ are known as  Jones projections for
$\B(I_1)\subset \A(I_1)$ and $\B(J_1)\subset \A(J_1)$ respectively.
They are related by an inner automorphism of $\B(I)$, which is the
following Lemma:

\blem\label{jones} Let $u=\lambda^{-1} E_I(w_1w_2^*) \in \B(I).$
Then $u$ is unitary and we have $e_1= ue_2u^*.$ \elem

\proof

First we have $w_1= u w_2$ and so  $e_1= ue_2u^*.$  Now compute
$$
uu^*= \lambda^{-2} E_I(w_1w_2^*) E_I(w_2w_1^*)=\lambda^{-2}
E_I(w_1w_2^* E_I(w_2w_1^*))=\lambda^{-1} E_I(w_1w_1^*)=1 \ ,
$$
where in the third equality we have used that $w_2^* E_I(w_2w_1^*)=\lambda
w_1^*$ and in the last equality that $ E_I(e_1)=\lambda$.

\qed

The following is proved in exactly the same way as in \cite{KLM}:
\begin{lemma}\label{index}
If $\B\subset \A$ is a   M\"{o}bius subnet such that $\mu_\A$ is
finite and $[\A:\B]<\infty.$ Then $\mu_\B =\mu_\A [\A:\B]^2.$
\end{lemma}

\section{Mutual information in the case of free fermions}\label{fermion}

\subsection{Basic representation of $ LU_r$ and free fermion net}

Let $H$ denote the Hilbert space $L^2(S^1; \mathbb C^r)$ of
square-summable $\mathbb{C}^r$-valued functions on the circle.  The
group $LU_r$ of smooth maps $S^1 \rightarrow U_r$, with $U_r$ the unitary group on $\mathbb C^r$, acts on $H$
multiplication operators.

Let us decompose $H= H_+ \oplus H_-$, where
$$
H_+ = \{\text{\rm functions whose negative Fourier coeffients vanish}\} \, .
$$
We denote by $p$ the Hardy projection from $H$ onto $H_+$.

Denote by $U_{\text{\rm res}}(H)$ the group consisting of unitary
operator $A$ on $H$ such that the commutator $[p, A]$ is a
Hilbert-Schmidt operator. Denote by $\text{\rm Diff}^+(S^1 )$ the
group of orientation preserving diffeomorphism of the circle. It
follows from Proposition 6.3.1 and Proposition 6.8.2 in \cite{PS}
that $LU_r $ and $\text{\rm Diff}^+(S^1 )$ are subgroups of
$U_{\text{\rm res}}(H)$. The basic representation of $ LU_r$ is the
representation
 on Fermionic Fock space $F_p = \Lambda(pH)
\otimes \Lambda ((1-p)H)^*$ as defined in \S 10.6 of \cite{PS}. For
more details, see \cite{PS} or \cite{W2}. Such a representation
gives rise to a graded net as follows. Denote by $\A_r (I)$ the von
Neumann algebra generated by $c(\xi )^\prime s$, with $\xi \in L^2
(I, \mathbb C ^r)$. Here $c(\xi ) = a(\xi ) + a(\xi )^*$ and $a(\xi )$
is the creation operator defined as in Chapter 1 of \cite{W2}. Let
$Z:F_p \rightarrow F_p $ be the Klein transformation given by
multiplication by 1 on even forms and by $i$ on odd forms. It
follows from \S 15 of chapter 2 of \cite{W2} that $\A_r$ is a graded
M\"{o}bius covariant net, and $\A_r$ will be called the {\it net of $r$
free fermions}. It follows from Prop. 1.3.2 of \cite{TL} that
$\A_r$ is strongly additive and \S 15 of chapter 2 of \cite{W2} that
$\mu_{\A_r}=1$.

Fix $I_i\in \PI, i=1,2$, and $I_1, I_2$ disjoint, that is $\bar
{I_1} \cap \bar {I_2} = \emptyset$, and $I=I_1\cup I_2$.

For bounded operators $A, B : F_p \rightarrow F_p$, we define $A^+
=\Gamma A \Gamma, \ A^- = A - A^+$, where $\Gamma$ is an operator on
$F_p$ given by multiplication by 1 on even forms and $-1$ on odd
forms. An operator $A$ is called even (resp. odd) if $A = A^+$
(resp. $A = A^-$).

We define a graded tensor product $\otimes_2$ by the following formula:
$$
A \otimes_2 B = A \otimes B^+ + A\Gamma \otimes B^- \ ,
$$
where $A \otimes_2 B$ is considered as an operator on Hilbert space
tensor product $F_p \otimes F_p$.

Let $A_1, A_2, B_1, B_2$ be even or odd operators, i.e. $\Gamma A_i
\Gamma = A_i$ or $-A_i, \ \Gamma B_i \Gamma = B_i$ or $-B_i, \ i =
1, 2$.  Define the degree $d(A) = 0$ or $1$ if $A$ is even or odd.

It follows from the definition of $\otimes_2$ that:
\begin{gather*}
(A_1 \otimes_2 B_1)^* = (-1)^{d(A_1)d(B_1)} A_1^* \otimes_2 B_1^* \ , \\
(A_1 \otimes_2 B_1) \cdot (A_2 \otimes_2 B_2) = (-1)^{d(B_1) d(A_2)}
A_1 A_2 \otimes_2 B_1 B_2 \ .
\end{gather*}
For $A \in \A_r (I_1), \ B \in \A_r(I_2)$, we define
$$
\omega(A \otimes_2 B) = \langle \Omega,AB \, \Omega \rangle
$$
where $\Omega$ is the vacuum vector in $F_p$.

\begin{lemma}\label{gradedtensor}
\item
{\rm(1)}  $\omega $ extends to a normal faithful state on the von
Neumann algebra $\{ A \otimes_2 B, \ A \in \A_r(I_1), \ B \in
\A_r(I_2) \}^{\prime\prime}$ (denoted by $\A_r(I_1) \hat \otimes_2
\A_r(I_2)$) on $F_p \otimes F_p$.  There exists a unitary operator
$U_1 : F_p \rightarrow F_p \otimes F_p$ such that:
$$
U_1 AB U_1^* = A \otimes_2 B \qquad \text{\rm for every} \qquad A \in
\A_r(I_1), \ B \in \A_r(I_2)\ .
$$
{\rm(2)} The unitary operator $U_1$ in (1)  can be chosen such that
$U_1^*(\Gamma \otimes \Gamma)U_1= \Gamma£¬$   hence
$U_1^*(\B(F_p)\otimes 1)U_1$ commutes with $Z\A_r(I_2)Z^{-1}$ and
therefore is $\Ad \Gamma$ invariant as a set and lies in
$\A_r(I_2')$ when $I_2$ is an interval.
\end{lemma}

\proof (1)  is proved in Prop. 2.3.1 of \cite{Xjw}. We note
that  by (2) the state $\omega_1\otimes_2 \omega_2$ defined in
Definition \ref{12state} is a normal state on type $\mathrm{III}$
factor $\A_r(I_1)\vee \A_r(I_2),$ and hence can be represented by a
unique vector $\psi$ in the positive cone associated with vector
state $\omega$ on $F_p.$ Since both $\omega_1\otimes_2 \omega_2 $
and $\omega$ are $\Ad \Gamma$ invariant, it follows that $\Ad
\Gamma$ preserves the positive cone, and $\omega_1\otimes_2 \omega_2
$ is also represented by $\Gamma\psi$. By uniqueness we must have
$\Gamma\psi=\psi.$ Now $U_1$ in (2) is uniquely fixed by the
condition $U_1\psi = \Omega\otimes \Omega,$ and it follows that
$U_1^* \Gamma \otimes \Gamma U_1= \Gamma£¬$  hence
$U_1^*\B(F_p)\otimes 1 U_1$ is $\Ad \Gamma$ invariant as a set,
graded commuting with $ \A_r (I_2)$ and therefore lies in
$\A_r(I_2')$ when $I_2$ is an interval by Cor. \ref{carp2}.

\qed

\brem\label{invf} If $\B\subset \A$ is a graded subnet, the proof of
(3) then applies to $\B$, and for any interval $I$, by  choosing
$I_{1n}\subset I_{2n}^c\subset I $ with $\cup_{n=1}^{\infty} I_{1n}
=I,$ we can get an increasing sequence of $\Ad \Gamma$ invariant (as
a set) finite dimensional type I factors $B_n$ such that $\cup_n
B_n$ is strongly dense in $\B(I)$.

\erem

\subsection{Mutual information for free fermions}
Let $I_1 , I_2\in\PI$ and $I= I_1\cap I_2$ as above.
\bdef\label{12state}
We set
$$\omega_1\otimes_2 \omega_2 (AB)= \langle  \Omega\otimes \Omega,
A\otimes_2\! B\, \Omega\otimes \Omega\rangle,\quad \forall A\in \A_r(I_1),\  B\in \A_r(I_2) \ .
$$
By (1) Lemma \ref{gradedtensor}  $\omega_1\otimes_2 \omega_2$
defines a normal state on $\A_r(I)$. We note that the restriction of $\omega_1\otimes_2
\omega_2$ to $\A_r(I_1)$ and $\A_r (I_2)$ is the same as
$\omega$.
\ede

The mutual information we will compute is
$S(\omega,\omega_1\otimes_2 \omega_2)$. When we wish to emphasize
the underlying net, we will also write the mutual information as
$S_{\A_r}(\omega,\omega_1\otimes_2 \omega_2).$ When $\B\subset \A_r$
is a subnet, we write $S_{\B}(\omega,\omega_1\otimes_2 \omega_2)$
the mutual information for the net $\B$ obtained by restricting
$\omega,\omega_1\otimes_2 \omega_2$ from $\A_r$ to $\B$. Note that
by (4) of Th. \ref{515} $S_{\B}(\omega,\omega_1\otimes_2
\omega_2)\leq S_{\A_r}(\omega,\omega_1\otimes_2 \omega_2)$.

$\omega$ on $\A_r(I)$ is  quasi-free state as studied by Araki in
\cite{arakif}. To describe this state, it is convenient to use
Cayley transform $V(x) = (x-i) / (x + i)$, which carries the (one point compactification of the) real
line onto the circle and the upper half plane onto the unit disk. It
induces a unitary map $$ {\displaystyle Uf(x)=\pi ^{-{\frac
{1}{2}}}(x+i)^{-1}f(V(x))}  $$ of $L^2(S^1, \mathbb{C}^r)$ onto
$L^2(\mathbb{R},\mathbb{C}^r )$. The operator $U$ carries the Hardy
space on the circle onto the Hardy space on the real line (cf.
\cite{RR}). We will use the Cayley transform to identify intervals on
the circle with one point removed to intervals on the real line.
Under the unitary transformation above, the Hardy projection on
$L^2(S^1, \mathbb{C}^r)$ is transformed to the Hardy projection on
$L^2(\mathbb{R},\mathbb{C}^r )$ given by £º
$$
Pf(x) = \frac{1}{2}f(x) + \int \frac{i}{2 \pi} \, \frac{1}{(x-y)}
f(y)dy \ ,
$$
where the singular integral is (proportional to) the Hilbert transform.

We write the kernel of the above integral transformation as $C$:
\begin{equation}
C(x,y)=\frac{1}{2} \delta(x-y) - \frac{i}{2 \pi} \,
\frac{1}{(x-y)}\, \ .
\end{equation}
The quasi free  state $\omega$ is determined by $$\omega
\big(a(f)^*a(g)\big)= \langle g, P f\rangle.$$ Slightly abusing our
notations, we will identify  $P$ with its kernel $C$ and simply
write  $$\omega \big(a(f)^*a(g) \big)= \langle g, C f\rangle.$$  $C$ will be
called   {\it covariance operator}.

Recall $I_i\in \PI, i=1,2$ , and $I_1, I_2$ are disjoint, that is
$\bar {I_1} \cap \bar {I_2} = \emptyset$, and $I=I_1\cup I_2$. We
assume that $I= (a_1,b_1)\cup (a_2,b_2)\cup ... \cup(a_n,b_n) $ in
increasing order.

\subsection{Computation of mutual information in finite dimensional case}\label{fd}

Choose finite dimensional subspaces $H_i$ of $L^2(I_i,\mathbb{C}_r),
i=1,2,$ and denote by $\Car(H_i)\subset \A(I_i)$ the corresponding
finite dimensional factors of dimensions $2^{2\dim H_i}$ generated
by $a(f), f\in H_i.$ Let $\rho_{12}$,  $\rho_1$,  $\rho_2$ be the density
matrices of the restriction of $\omega$ to $\Car(H_1)\otimes_2
\Car(H_2)$, $\Car(H_1)$, $\Car(H_2)$ respectively, and
$\rho_1\otimes_2\rho_2$ of the restriction of
$\omega_{1}\otimes_2\omega_2$ to $\Car(H_1)\otimes_2 \Car(H_2)$. Our
goal in this section is to compute the relative entropy
$S(\rho_{12},\rho_1\otimes_2\rho_2).$

Note that since  $\Car(H_1)$ is type I factor, $\Ad \Gamma$ acts on
 $\Car(H_1)$ by an inner automorphism $\Ad u, u\in \Car(H_1).$ Since
 $\Ad u$ has order two, by suitably choosing phase factor we can
 assume that $u^2=1.$ Note that $\Gamma u \Gamma = u^3= u,$ so $u$ is
 even, and $\Gamma u$ commutes with  $\Car(H_1).$ So $\Gamma u
 \otimes B^-, 1\otimes B^+$ generates a type I factor $\widetilde{\Car}(H_2)$  isomorphic to
 $\Car(H_2),$  and commuting with $\Car(H_1)\otimes 1.$ It follows
 that  $\Car(H_1)\otimes_2 \Car(H_2)= \Car(H_1)\otimes
 \widetilde{\Car}(H_2)$.  Let us show that
 $\omega_{1}\otimes_2\omega_2,$ when restricting to  $\Car(H_1)\otimes
 \widetilde{\Car}(H_2)$, is the tensor product state $\rho_1\otimes
 \rho_2',$ where $\rho_1, \rho_2'$ denote the restriction of $\omega$ to
 $\Car(H_1)$, $\widetilde{\Car}(H_2)$ respectively. Since $\omega_{1}\otimes_2\omega_2$
 clearly agrees with $\rho_1\otimes
 \rho_2'$  on $A\otimes B^+$,   it is sufficient to check that
 $$
 \omega_{1}\otimes_2\omega_2(a \tilde b)= \omega(a)\omega(\tilde b),\ \ \forall a\in
 \Car(H_1)\otimes 1,\ \tilde b = \Gamma u \otimes b^-,\ b^-\in
 \Car(H_1)^-\ .
 $$
The left-hand side of the above is
$$
\langle \Omega, au\Omega\rangle\langle \Omega,b^- \Omega\rangle=0
$$
and the right-hand side is
$$
\langle \Omega, a\Omega\rangle\langle \Omega,ub^- \Omega\rangle=0
$$
since $u$ is even. We also note that $\omega$ restricted to
$\widetilde{\Car}(H_2)$ is the same as  $\omega$ restricted to
$\Car(H_2)$ under the natural isomorphism of $\widetilde{\Car}(H_2)$
with $\Car(H_2)$.

So we have shown the analog of \eqref{remu} in this graded local context:

\bprop\label{finitev}
$$S(\rho_{12},\rho_1\otimes_2\rho_2)
=S(\rho_{1})+S(\rho_{2}) -S(\rho_{12}) \ .
$$
\eprop
Now we turn to the computation of von Neumann entropy $S(\rho_{1}).$
Let $p_1$ be the projection  onto the finite dimensional subspace
$H_1$ of $L^2(I_1,\mathbb{C}_r)$. $\rho_1$ on $\Car(H_1)$ is quasi
free state given by covariance operator $C_{p_1}=p_1Cp_1.$ Let $K$
be the operator such that
$$
(1+\exp(-K))= C_{p_1}
$$
Since $K$ is self adjoint, we can  choose an orthonormal basis
$\psi_i , \ 1\leq i\leq \dim H_1$ of $H_1$  such that $K\psi_i =
\lambda_i \psi,$ where $\lambda_i$ are real eigenvalues of $K$.

$\Car(H_1)$ acts on the Fermionic Fock space $F(H_1)$.
Let $$ K_1:=
\sum_{i}  \lambda_i a(\psi_i)^* a(\psi_i)\ .$$
According to \cite{arakif} and \cite{Casf}, the density matrix of
$\rho_1$  (still denoted by $\rho_1$) as an operator on $F(H_1)$ is
given by the following
$$
\rho_1=  c \exp(-K_1) \ ,
$$
where $c^{-1}= \Tr \big(\exp(-K_1)\big)$.

By a simple computation we find that $\Tr
\big(\exp(-K_1)\big)=\textrm{det}(1+e^{-K})$ and

\ben\label{srho}
 S(\rho_1)= \Tr (\rho_1\ln \rho_1)=
\Tr\big((1-C_{p_1})\,\log (1-C_{p_1})+C_{p_1} \,\log C_{p_1}\big) \ . \een

\bdef\label{sigma}

Let $\mathbf{P}_i$ be projections from $L^2(I,\mathbb{C}^r)$ onto
$L^2({I_i},\mathbb{C}^r),$ and $C_i= \mathbf{P}_iC \mathbf{P}_i,
i=1,2$.

Let
\begin{multline*}
\sigma_C = \mathbf{P_1}\big (C\ln C +(1-C)\ln (1-C)\big)\mathbf{P_1}- \big(C_1
\ln C_1 +  (\mathbf{P_1}-C_1)\ln( \mathbf{P_1}-C_1)\big)+ \\
\mathbf{P_2}\big(C\ln C +(1-C)\ln (1-C)\big)\mathbf{P_2}- \big(C_2 \ln
C_2+(\mathbf{P_2}-C_2)\ln( \mathbf{P_2}-C_2)\big)
\end{multline*}
 and $\sigma_{C_p}$ be
the same as in the definition of $\sigma_C $ with $C$ replaced by
$C_p= pCp,$ if $p$ is a projection commuting with $\mathbf{P}_1.$

\ede

Denote by $p$ the projection from $L^2(I,\mathbb{C}^r)$ onto
$H_1\oplus H_2$. By Prop. \ref{finitev} and equation (\ref{srho}) we
have proved the following

\bprop\label{finiter}

$$S(\rho_{12},\rho_1\otimes_2\rho_2)
=\Tr (\sigma_{C_p}) \ .
$$
\eprop

\subsection{Inequality from operator convexity}

The proof of the following result can be found in \cite{Cal} (See
Th. 2.6 and Th. 4.19 of \cite{Cal}):

\bthm \label{SD}(1) For all operator convex functions $f$ on $\mathbb R$, and all
orthogonal projections $p$, we have $p f(pAp)p \leq p f(A)p$ for every selfadjoint operator $A$; (2) $f(t) = t
\ln(t)$ is operator convex. \ethm

(1) of the above Theorem is known as Sherman-Davis Inequality.  It
in instructive to review the idea of the proof of (1) which is also
used in the proof of Th. \ref{AB}: Consider the selfadjoint unitary operator $U^p =2p-I$;  by operator convexity we have
$$
f\big(\frac{1}{2} A +\frac{1}{2} U^p A U^p\big) \leq \frac{1}{2} f(A)
+\frac{1}{2} f(U^p A U^p) \ .
$$
Now notice that $$\frac{1}{2} A +\frac{1}{2} U^p A U^p = A_p +
A_{1-p},\quad f(U^p A U^p)= U^p f(A) U^p\ ,$$
where $A_p = pAp$,
and the inequality follows.

For (2), see e.g. \cite{Cal}.

\blem\label{ine1} (1) $$ S(\omega, \omega_1\otimes_2 \omega_2)=
\lim_{p\rightarrow 1} \Tr(\sigma_{C_p})\geq \Tr(\sigma_C)
$$ where $p\rightarrow 1$ strongly. \par
(2) The mutual information for $r$ free fermion net  is $r$ times
the mutual information for $1$ free fermion net. \elem

\proof

(1): The first follows from Prop. \ref{finiter} and (2) of Th.
\ref{515}. To prove the inequality, we use the fact that $x\ln x$ is
operator convex, and so $\mathbf{P_1} C\ln C\mathbf{P_1}\geq C_1 \ln
C_1$, and similarly with $C$ replaced by $1-C$ by Th. \ref{SD}. It
follows that $\sigma\geq 0, \sigma_p\geq 0.$ Since $\sigma_p$ goes
to $\sigma$ strongly as $p\rightarrow 1$ strongly, the inequality
follows.\par

(2): For the case of $r$ free fermions, the trace in (1) is over
$L^2(\mathbb{R}, \mathbb{C}^r)$  which is $r$ direct sum of the
Hilbert space $L^2(\mathbb{R}, \mathbb{C}),$ and (2) follows.

\qed

We shall prove later that the inequality in the above Lemma is
actually an equality. It would follow if one can show that
$\sigma_{C_p}$ goes to $\sigma_C$ in tracial norm. This is not so
easy, and we note that  $\mathbf{P_1} \big(C\ln C +(1-C)\ln
(1-C)\big)\mathbf{P_1}$ is not trace class. To overcome this difficulty
and to compute the mutual information we prove the reverse
inequality by applying Lieb's joint convexity and regularized kernel
as in the next two sections.

\subsection{Reversed inequality from Lieb's joint convexity}

We begin with the following Lieb's Concavity Theorem:

\bthm \label{liebc} (1) For all $m \times  n$ matrices $K$, and all
$0\leq t \leq  1$, the real valued map  given by $(A, B) \rightarrow
\Tr(K^* A^{1-t} KB )$ is concave where $A,B$ are non-negative
$m\times m$ and $n\times n$ matrices respectively;
\par

(2) If $A\geq 0, B\geq 0$ and $K$ is trace class,  then
$$(A, B) \rightarrow \Tr(K^* A^{1-t} KB), \quad 0\leq t\leq 1,$$
is jointly concave;

(3) If  $A\geq \epsilon I, B\geq \epsilon I,\epsilon>0$ and $K$ is
trace class, then
$$(A, B) \rightarrow \Tr(K^* A\ln A K- K^* AK \ln B) $$
is jointly convex; \ethm

\proof

(1) is proved in Th. 6.1 of \cite{Cal}. (2) follows from (2) by
functional calculus. To prove (3), we note that
$$
\Tr(K^* A\ln A K - K^* AK \ln B)= \lim_{t\rightarrow 0} \frac{ \Tr(K^*
A^{1-t} KB) - \Tr (K^*AK)}{t-1}
$$ and (3) follows from (2).

\qed

\blem\label{trace1} Assume that $S$ is trace class, then $\Tr
(ST)=\Tr (TS)$ where $T$ is any bounded operator, and if the sequence of bounded operators
$T_n\rightarrow T$ strongly, then $\Tr (ST_n)\rightarrow \Tr (ST) $.
\elem

\proof The equality is proved in \cite{RS}. Let $e_i$ be  an
orthonormal basis, and $S=U|S|$ be the polar decomposition of $S$.
Then
$$
\Tr (T_nS)= \sum_{i} \langle e_i, T_n U|S|^{1/2} |S|^{1/2}e_i\rangle \ .
$$
Note that
$$
|\langle e_i,T_n U|S|^{1/2} |S|^{1/2}e_i\rangle|\leq ||T_n
U|S|^{1/2} e_i||\, || |S|^{1/2}e_i||\leq c \langle e_i,|S|e_i\rangle,\
\forall i \ ,
$$
where $c$ is a constant, so the last part of the Lemma follows
by Lebesgue dominated convergence theorem.

\qed

\blem\label{ln} Suppose that $K\geq \epsilon I,L\geq \epsilon
I,\epsilon
>0 $ and $K-L$ is trace class. Then $\ln K- \ln L$ is trace class.
\elem \proof Note that
$$
\ln K=-\int_{0}^\infty \left(\frac{1}{K+t}-\frac{1}{1+t}\right)dt,\quad \ln
L=-\int_{0}^\infty \left(\frac{1}{L+t}-\frac{1}{1+t}\right)dt\ .
$$
Hence
$$
\ln K- \ln L= -\int_{0}^\infty \left(\frac{1}{K+t}-\frac{1}{L+t}\right)dt
=\int_{0}^\infty \left(\frac{1}{L+t} ( K-L) \frac{1}{K+t}\right)dt \ .
$$
We have
\begin{align*}
||\ln K- \ln L||_1\leq \int_{0}^\infty \left\Vert\frac{1}{L+t} ( K-L)
\frac{1}{K+t}\right\Vert_1 dt \\
\leq \int_{0}^\infty \left\Vert\frac{1}{L+t}\right\Vert \left\Vert
K-L\right\Vert_1  \left\Vert \frac{1}{K+t}\right\Vert  dt \leq ||K-L||_1\, \epsilon ^{-1}\ ,
\end{align*}
where $||\cdot||_1$ denotes tracial norm.

 \qed

\bthm\label{AB} Let $A\geq \epsilon ,\epsilon>0, B:=
\mathbf{P_1}A\mathbf{P_1}+ \mathbf{P_2}A\mathbf{P_2},$ where
$\mathbf{P_1}$ is a projection, $\mathbf{P_1}+\mathbf{P_2}=1$,  and
$p$ is a finite rank projection commuting with $\mathbf{P_1}.$
Assume that $A-B$ is trace class. Then
$$
\Tr \big(A(\ln A- \ln B)\big) \geq \Tr \big(A_p (\ln A_p- \ln B_p)\big) \ .
$$
\ethm

\proof

Apply Th. \ref{liebc} to $A, B$ and unitary $U^p = 2P- I,$ with
$f(A,B,K)= \Tr (K^* A\ln A K- K^* AK \ln B),$ $K$ is a finite rank
projection, we have
$$
f\left(\frac{1}{2}(A+ U^p A U^p), \frac{1}{2}(B+ U^p B U^p), K\right) \leq
\frac{1}{2} f(A,B,K) + \frac{1}{2} f(U^p A U^p, U^p B U^p,K) \ .
$$
Note that
$$
f\left(\frac{1}{2}(A+ U^p A U^p), \frac{1}{2}(B+ U^p B U^p),
K\right)=f(A_p+A_{1-p}, B_p+B_{1-p}, K)
$$
and

$$
f(A_p+A_{1-p}, B_p+B_{1-p}, K)=
 \Tr \big(K (A_p \ln A_p + A_{1-p} \ln A_{1-p})K- K(A_p+A_{1-p})K \ln
(B_p +B_{1-p})\big)
$$
and

\begin{multline*}
\frac{1}{2}f(A,B,K) + \frac{1}{2} f(U^p A U^p, U^p B U^p,K)= \\
\frac{1}{2} \Tr (K A\ln A\, K- K AK \ln B) + \frac{1}{2} \Tr (K U^p
A\ln A\, U^p K - KU^p A U^p K U^p \ln B\, U^p)\ .
\end{multline*}
Observe that
$KU^p = [K,2p]+ U^p K$ and  $K\ln B= [K,\ln B] + \ln B\,
K, K\ln (B_p +B_{1-p})= [K,\ln (B_p +B_{1-p})] + \ln (B_p
+B_{1-p})K.$ We will let $K\rightarrow 1$ strongly eventually.  Up
to terms that go to $0$ as $K\rightarrow 1$ strongly, we can freely
permute $K$ and $U^p$. By permuting $K$ with $\ln(B_p +B_{1-p})$ and
$\ln B$ on the left hand side and righthand side of the above
inequality respectively, we get terms on the left hand side of the
above inequality
$$
-\Tr \big(K (A_p+A_{1-p}) [K, \ln (B_p+ B_{1-p})]\big)
$$
and on the right hand side of the above inequality
$$
- \frac{1}{2}\Tr \big(K A [K, \ln B] +K U^p A [K,\ln B] U^p\big) = -\Tr \big(K
(pA[K,\ln B] p +(1-p)A[K,\ln B](1-p)\big).
$$
Up to terms that go to $0$ as  $K\rightarrow 1$ strongly, we have that
$$
-\Tr \Big(K \big(pA[K,\ln B] p +(1-p)A[K,\ln B](1-p)\big)\Big)
$$
is equal to
$$
-\Tr \Big(K \big(pAp[K,\ln B] +(1-p)A(1-p)[K,\ln B](1-p)\big)\Big) \ .
$$
This is the same as $$ -\Tr \big(K (A_p+A_{1-p}) [K, \ln (B_p+ B_{1-p})]\big)
$$ up to  terms that go to $0$ as  $K\rightarrow 1$ strongly since
$\ln B- \ln (B_p + B_{1-p})$ is trace class by Lemma \ref{ln}, and
$p$ is of finite rank.

So by permuting $K$ with $U^p$ and $\ln B$, $\ln(B_p +B_{1-p})$ in
the inequalities above, by  Lemma \ref{trace1}, use $\ln B- \ln (B_p
+ B_{1-p})$ is trace class by Lemma \ref{ln}, and $p$ is finite
rank, we have that up to terms which go to zero as $K\rightarrow 1$
strongly,
\begin{multline*}
\Tr \Big(K \big(A_p (\ln A_p -\ln B_p) + A_{1-p}(\ln A_{1-p} - \ln B_{1-p}
)\big)\Big) \\
\leq \Tr \Big(K \big(pA (\ln A-\ln B) p + (1-p) A(\ln A - \ln B
)(1-p)\big)\Big)\ .
\end{multline*}
By Lemma \ref{ln}, $\ln A_p -\ln B_p, \ln A_{1-p} - \ln B_{1-p}, \ln
A - \ln B$ are trace class operators, and by Lemma \ref{trace1}
 let $K\rightarrow 1$ strongly, we get
$$
 \Tr \big(A_p (\ln A_p -\ln B_p)\big) + \Tr \big(A_{1-p }(\ln A_{1-p} -\ln
 B_{1-p})\big)
 \leq \Tr A(\ln A -\ln B) \ .
$$
As in the proof of Lemma \ref{ine1}, by operator convexity of $xlnx$
we have $$(1-p)A_{1-p }\ln A_{1-p}(1-p) \geq B_{1-p}\ln B_{1-p}$$
and so
$$
\Tr \big(A_{1-p }(\ln A_{1-p} -\ln B_{1-p}))= \Tr (A_{1-p }\ln A_{1-p}-
B_{1-p}\ln B_{1-p}\big)\geq 0
$$
and the theorem is proved. \qed

\subsection{Regularized Kernel for one free fermion case}
Note that, by Lemma \ref{ine1}, the mutual information for $r$ free
fermion net  is $r$ times the mutual information for $1$ free
fermion net. In this section we will determine the mutual
information for $1$ free fermion net.

This section is inspired by formal computations in \cite{Casf}. The
regularization is also  motivated by Th. \ref{AB} which  applies to
strictly positive operators.

Recall that the Hardy projection on $L^2(\mathbb{R}, \mathbb{C})$ is
given by £º
$$
Pf(x) = \frac{1}{2}f(x) + \int \frac{i}{2 \pi} \frac{1}{(x-y)}
f(y)dy \ ,
$$
where the integral is the singular integral or Hilbert transform.

We write the kernel of the above integral transformation as $C$.
\begin{equation}
C(x,y)=\frac{1}{2} \delta(x-y) - \frac{i}{2 \pi} \,
\frac{1}{(x-y)}\ \ .
\end{equation}
Recall $I_i\in \PI, i=1,2$ , and $I_1, I_2$ are disjoint, that is
$\bar {I_1} \cap \bar {I_2} = \emptyset$, and $I=I_1\cup I_2$. We
assume that $I= (a_1,b_1)\cup (a_1,b_1)\cup ... \cup(a_n,b_n) $ in
increasing order.

Then resolvent of $C$ as restriction of an operator on
$L^2(I,\mathbb{C})$
   \begin{equation}
   R^0(\beta)=(C-1/2+\beta)^{-1}\label{nu}\equiv\left( - \frac{i}{2\pi}\frac{1}{x-y} +\beta \,\delta(x-y)\right)^{-1}
   \end{equation}
 has the following expression  (\cite{singu} or Page 133 of
 \cite{MSG}):
\begin{equation}
R^0(\beta) =\left(\beta^2-1/4 \right)^{-1}
\left(\beta\,\delta(x-y)\, +\frac{i }{2\pi} \frac{e^{-\frac{i}{2\pi}
\log\left(\frac{\beta-1/2}{\beta+1/2}\right)\, (Z(x)-Z(y)) }}{x-y}
\right)\,,\label{reso}
\end{equation}
where
\begin{equation}
Z(x)=\log\left(-\frac{\prod_{i=1}^n (x-a_i)}{\prod_{i=1}^n
(x-b_i)}\right) \,. \label{bfbf}
\end{equation}
It is useful to consider the following regularized operator: Let
$\epsilon_0>0$, and $E:= \frac{C+\epsilon_0}{1+2\epsilon_0}.$ Note
that
$$\frac{1+\epsilon_0}{1+2\epsilon_0}\geq E\geq
\frac{\epsilon_0}{1+2\epsilon_0}$$
Then we have
$$
E\ln E+ (1-E)\ln (1-E) = \int_{\frac{1}{2}}^\infty
\left[\big(\beta-\frac{1}{2}\big) \big(R_E(\beta)-R_E(-\beta)\big)-
\frac{2\beta}{\beta+\frac{1}{2}}\right] d\beta \ ,
$$
where $R_E(\beta)= \frac{1}{E-\frac{1}{2}+\beta}$\ .

We note that the integral above is absolutely convergent in norm.
This can be seen as follows: the integrand is
$$
\big(\beta-\frac{1}{2}\big) \big(R_E(\beta)-R_E(-\beta)\big)-
\frac{2\beta}{\beta+\frac{1}{2}}=\frac{\beta/2-
2\beta(E-\frac{1}{2})^2}{[(E-\frac{1}{2})^2-\beta^2](\beta+\frac{1}{2})} \ .
$$
For $1/2\leq \beta\leq 1$, since
$$\frac{1+\epsilon_0}{1+2\epsilon_0}\geq E\geq
\frac{\epsilon_0}{1+2\epsilon_0}\ ,$$ we have
$$
\left\Vert\frac{2\beta}{\beta+\frac{1}{2}}=\frac{\beta/2-
2\beta(E-\frac{1}{2})^2}{[(E-\frac{1}{2})^2-\beta^2](\beta+\frac{1}{2})}\right\Vert\leq
\frac{3\beta}{\beta+ 1/2}\left(\frac{\epsilon_0}{1+2\epsilon_0}\right)^{-2} \ .
$$
On the other hand
$$\left\Vert\frac{2\beta}{\beta+\frac{1}{2}}=\frac{\beta/2-
2\beta(E-\frac{1}{2})^2}{[(E-\frac{1}{2})^2-\beta^2](\beta+\frac{1}{2})}\right\Vert$$
is bounded by $\frac{1}{\beta^2}$ when $\beta$ is large.

To evaluate the above integral using resolvent, let
$t=\beta(1+2\epsilon_0)$ we get
\begin{multline*}
E\ln E+ (1-E)\ln (1-E)\\=\int_{\frac{1}{2}(1+2\epsilon_0)}^\infty
\left[\left(\frac{t}{1+2\epsilon_0}-\frac{1}{2}\right) \big(R(t)-R(-t)\big)-
\frac{2t}{t+\frac{1}{2}(1+2\epsilon_0)} \frac{1}{1+2\epsilon_0}\right]
d\beta.
\end{multline*}
Now we determine the kernel $K_1^{\epsilon_0}(x,y), x,y\in I_1$ of
$$
\mathbf{P}_1 E\ln E+ (1-E)\ln (1-E)\mathbf{P}_1-E_1\ln E_1-
(\mathbf{P}_1-E_1)\ln (\mathbf{P}_1-E_1) \ .
$$

\blem\label{gc} Suppose $f\in C^1(I_1\times I_1)$ and $f(x,y)=
-f(y,x)$. Let $g(x,y)= \frac{f(x,y)}{x-y} $ if $x\neq y$ and
$g(x,x)= \frac{\partial f}{\partial x}. $ Then $g(x,y)$ is
continuous on $I_1\times I_1.$ \elem \proof

It is enough to check continuity at $(y,y), y\in I_1.$ Since  $f\in
C^1(I_1\times I_1),$ we can write $f(x',y')=\frac{\partial
f}{\partial x}(y,y)(x'-x) +\frac{\partial f}{\partial y}(y,y)(y'-y)
+o(x'-y')=\frac{\partial f}{\partial x}(y,y)(x'-y')  +o(x'-y') $
where in the second $=$ we have used  $f(x,y)= -f(y,x)$ and hence
$\frac{\partial f}{\partial x}(y,y)= -\frac{\partial f}{\partial
y}(y,y)$. It follows that $\lim_{(x',y')\rightarrow (y,y)} g(x',y')=
g(y,y).$

\qed

We shall denote by $Z_{I,I_1}(x) = Z_{I}(x)- Z_{I_1}(x).$ Even
though both $Z_{I}(x )$ and $Z_{I_1}(x)$ are singular when $x$ is
close to the boundary of its domain, it is crucial that
$Z_{I,I_1}(x)$ is a smooth function on the closure of $\bar I_1.$

\blem\label{Gc} Let
$$G(t,x,y)=
\frac{\sin\left(\frac{1}{2\pi}\ln\left(\frac{t-\frac{1}{2}}{t+\frac{1}{2}}\right)\big(Z_I(x)-Z_I(y)\big)\right)-
\sin\left(\frac{1}{2\pi}\ln\left(\frac{t-\frac{1}{2}}{t+\frac{1}{2}}\right)\big(Z_{I_1}(x)-Z_{I_1}(y)\big)\right)}{x-y}
$$ if $x\neq y$ and $G(t,x,x)= \frac{1}{2\pi}\ln
\left(\frac{t-\frac{1}{2}}{t+\frac{1}{2}}\right) \big(Z_I'(x)- Z_{I_1}'(x)\big)$,
$t>\frac{1}{2}$.

Then $G(t,x,y)$ is continuous on $(\frac{1}{2},
\infty)\times I_1\times I_1$ and
$$
|G(t,x,y)|\leq
\left|\frac{1}{2\pi}\ln\left(\frac{t-\frac{1}{2}}{t+\frac{1}{2}}\right)\right| M, \quad
(t,x,y)\in (\frac{1}{2}, \infty)\times I_1\times I_1\ ,
$$ where $M$ is a constant.

\elem

\proof The continuity of $G$ follows from Lemma \ref{gc}. To prove
the inequality, we note that
$$
|G(t,x,y)|\leq
\left|\frac{1}{2\pi}\ln\left(\frac{t-\frac{1}{2}}{t+\frac{1}{2}}\right)\right|\,
\left|\frac{Z_{I,I_1}(x)-Z_{I,I_1}(y)}{x-y}\right| \ .
$$
We note that $Z_{I,I_1}(x)-Z_{I,I_1}(y)$ is smooth on $\bar
I_1\times \bar I_1,$ and apply Lemma \ref{gc} we have proved the
inequality.

\qed

By Lemma \ref{Gc}, we have that the kernel before Lemma \ref{gc} is given
by
$$
K^{\epsilon_0}_1(x,y) = \frac{-1}{\pi}
\int_{\frac{1}{2}(1+2\epsilon_0)}^\infty
\frac{\big(\frac{t}{1+2\epsilon_0}-\frac{1}{2}\big)}{t^2-1/4} G(t,x,y)dt \ .
$$

\blem\label{Ku} (1) $K^{\epsilon_0}_1(x,y)$ is continuous, uniformly
bounded and converges uniformly on $I_1\times I_1$ to $K^{0}_1(x,y)$
as $\epsilon_0$ goes to $0$; \par

(2) The kernel of $$\mathbf{P}_1C\ln C+ (1-C)\ln
(1-C)\mathbf{P}_1-C_1\ln C_1- (\mathbf{P_1}-C_1)\ln
(\mathbf{P}_1-C_1)$$ is given by the bounded continuous function
$K^{0}_1(x,y),$ and moreover its trace is given by
$$
\int_{I_1} K^{0}_1(x,x) dx= \lim_{\epsilon_0\rightarrow 0}\int_{I_1}
K^{\epsilon_0}_1(x,x) dx ;$$\par

(3)
$$\int_{I_1} K^{0}_1(x,x) dx = \frac{1}{12}\sum_{
(a_i,b_i)\in I_2, (a_j,b_j)\in I_1}
\ln\left(\frac{(a_j-a_i)(b_j-b_i)}{(b_j-a_i)(a_j-b_i)}\right) \ .$$

\elem

\proof (1): It is clear  $K^{\epsilon_0}_1(x,y)$ is continuous
and uniformly bounded by Lemma \ref{Gc}. By Lemma \ref{Gc} again

$$
|K^{\epsilon_0}_1(x,y)-K^{0}_1(x,y)|\leq \frac{1}{2\pi^2} M
\int_{1/2}^\infty
\left|\frac{\big(\frac{t}{1+2\epsilon_0}-\frac{1}{2}\big)}{t^2-1/4}\chi_{\left(\frac{1}{2}(1+2\epsilon_0),\infty\right)}
-\frac{1}{t+1/2}\right| \left| \ln\left(\frac{t-\frac{1}{2}}{t+\frac{1}{2}}\right)\right| dt
$$
where $\chi_{\left(\frac{1}{2}(1+2\epsilon_0),\infty\right)}$ denotes the
characteristic function. We note that the integrand above is bounded
and when $t$ is large decays like a constant multiply by
$\frac{1}{t^2}$.

The uniform convergence now follows by Lebesgue's  dominated
convergence theorem. \par

(2): Note that as $\epsilon_0$ goes to $0$, $\mathbf{P}_1 E\ln E+
(1-E)\ln (1-E)\mathbf{P}_1-E_1\ln E_1-(\mathbf{P}_1-E_1)\ln
(\mathbf{P}_1-E_1)$ converges to
$$\mathbf{P}_1C\ln C+ (1-C)\ln (1-C)\mathbf{P}_1-C_1\ln C_1-
(\mathbf{P}_1-C_1)\ln (\mathbf{P}_1-C_1)$$ strongly. (2) now follows
from (1) and \cite{B} which contains more general results on the trace
of operators with integrable kernels.
\par

(3): By Lemma \ref{Gc} and (2) we have
\begin{multline*}
\int_{I_1} K^{0}_1(x,x) dx = \\ \frac{-1}{2\pi^2}\int_{\frac{1}{2}}^1 \frac{1}{t+1/2} \ln\left(\frac{t-\frac{1}{2}}{t+\frac{1}{2}}\right) \left(\sum_{
(a_i,b_i)\in I_2, (a_j,b_j)\in I_1}
\ln\left(\frac{(a_j-a_i)(b_j-b_i)}{(b_j-a_i)(a_j-b_i)}\right)\right)dt \ .
\end{multline*}
To finish the proof we just need to show
$\frac{-1}{2\pi^2}\int_{\frac{1}{2}}^1 \frac{1}{t+1/2}
\ln\left(\frac{t-\frac{1}{2}}{t+\frac{1}{2}}\right)= 1/12.$ By change of
integration variable to $u=\ln\left(\frac{t-\frac{1}{2}}{t+\frac{1}{2}}\right)$
it is sufficient to check that
$$
\int_{-\infty}^0 \frac{ue^u}{1-e^u} du = \frac{-1}{6\pi^2} \ .
$$
Since the anti-derivative of $\frac{ue^u}{1-e^u}$ is $-\textbf{Li}_2
(e^u) - u\ln (1-e^u)$ where $\textbf{Li}_2(x):= \sum_{k=1}^{\infty}
\frac{x^k}{k^2}$ is the dilogarithm, the desired equality follows
from
$$
\sum_{k=1}^\infty \frac{1}{k^2}=\frac{\pi^2}{6}\ .
$$
\qed
\brem\label{i2} We note that the previous Lemma works in exactly the
same way when we replace $I_1$ by $I_2$, and $\mathbf{P}_1$ by
$\mathbf{P}_2$. \erem

\subsection{The proof of Theorem \ref{main1}}

\bdef\label{G}
 If $I=(a_1, b_1)\cup (a_2, b_2)\cup ... \cup
(a_n, b_n)$ in increasing order,  define

$$G(I):=\frac{1}{6} \left(
\sum_{i,j}\log|b_i-a_j|-\sum_{i<j} \log|a_i-a_j| -\sum_{i<j}
\log|b_i-b_j| \right) \ .
$$
\ede

\bthm\label{main1} Let $I=(a_1, b_1)\cup (a_2, b_2)\cup ... \cup
(a_n, b_n)\in \PI$ and $I_1\cup I_2 =I, \bar I_1\cap \bar
I_2=\emptyset$.
Then
$$ S_{\A_r}(\omega, \omega_1\otimes_2\omega_2)
=r\big(G(I_1)+G(I_2)-G(I_1\cup I_2)\big)\ .
$$

\ethm

\proof

By Lemma \ref{ine1} it is sufficient to prove $r=1$ case.

Recall that $E:= \frac{C+\epsilon_0}{1+2\epsilon_0}.$
 Apply Theorem \ref{AB} to $A= E$ and
$A=(1-E)$ respectively, we have
$$
\Tr\sigma_{E}\geq \Tr\sigma_{E_p} \ .
$$
Now let $\epsilon_0$ go to $0$ and by (2), (3) of Lemma \ref{Ku},
Lemma \ref{ine1} and Remark \ref{i2},   Theorem \ref{main1} is
proved.

\qed

\section{Subnets of free fermion nets and their finite index extensions}\label{descent}

\subsection{Formal properties of entropy for free fermion nets  and their subnets}
In the previous section we use Cayley transformation to identify
punctured circle with real line as a  tool to compute relative
entropy. Now we return to general discussion about formal properties
of entropy, and it is now convenient to be back to intervals on the
circle.  Let $I\in \PI$ be disjoint union of intervals on the
circle. Explicitly we write $I= (a_1,b_1)\cup (a_2,b_2)\cup ...
\cup(a_n,b_n) $ in anti-clockwise order on the unit circle. We note
that relative entropies as computed in Th. \ref{main1} is invariant
under  $\Mob$ transformations on the circle. The results of this
section are inspired by \cite{Casc}.\par

By  Theorem \ref{main1}, we have $F_{\A_r}(A,B) := S(\omega,
\omega_A \otimes_2 \omega_B)< \infty$ where $A, B$ are union of
disjoint intervals . When no confusion arises, we will simply write
$F_{\A_r}(A,B)$ as $F(A,B)$.

We can extend the definition mutual information to more general
union of disjoint intervals by the following
$$
F(A\cup B, A\cup C)= F(A, B\cup C) + F(B, C) - F(A,C)- F(A,B) \ .
$$

\bthm\label{soft1}

(1) $$ F(A\cup B, A\cup C) \geq 0;$$ $F(A\cup B, A\cup C)$ is
continues from inside;\par

(2)\begin{multline*} F(A,B) + F(A,C) + F(A \cup B,A \cup C) + F(A \cap C,A \cap B) \\ =
F(B,C) + F(A,B \cup C) + F(A,B \cap C).
\end{multline*}

(3) There exists function $G: \PI \rightarrow \mathbb{R}$ such that
$$
F(A,B) = G(A)+G(B)- G(A\cup B) -G(A\cap B) \ .
$$
Such $G$ is uniquely determined by its value on connected open
intervals;

(4) One can choose $G(a,b)=\frac{r}{6} \ln |b-a|$ in (3) for the $r$
free fermion net $\A_r$, and such a choice determines
$$G(I) =
\frac{r}{6} \left( \sum_{i,j} \ln |b_i-a_j| - \sum_{i<j} \ln |a_i-a_j|
-\sum_{i<j} \ln |b_i-b_j|\right)
$$
for $I= (a_1,b_1)\cup (a_2,b_2)\cup ...
\cup (a_n,b_n)$ on unit circle with anti-clockwise order;
\par

(5) $ F(A\cup B, A\cup C)= F(A\cup B,C)- F(A,C) =  F(B,A \cup C)-
F(B,A);$ In particular  $F(A\cup B, A\cup C)$ increases with
$B,C;$\par

(6) If $\B\subset \A$ is a graded  subnet, then (1), (2), (3) is also
true for the system of mutual information associated with $\B$.\ethm

\proof

(1) and (5) for free fermions can be checked by using explicit
formulas in Th. \ref{main1}, but here we present general arguments
which will also works for other cases such as subnets of  free
fermions.

Choose increasing sequence of finite dimensional factors $I_{A_n},
I_{B_n},$ invariant under the conjugate action of $\Gamma$ such that
$(\bigcup_n I_{A_n})''= \A_r(A)$, $(\bigcup_n I_{B_n})''= \A_r(B)$,
and denote by $\rho_{A_nB_n}, \rho_{A_n}\otimes_2 \rho_{B_n}$ the
restrictions of $\omega$ and $\omega_1\otimes_2\omega_2$ to
$I_{A_n}\vee I_{B_n}$ respectively. Let $\rho_{A_n}$ and
$\rho_{B_n}$ be  the restrictions of $\omega$ to $I_{A_n}$and $
I_{B_n}$ respectively.

By  Prop. \ref{finiter}
$$
S(\rho_{A_nB_n}, \rho_{A_n}\otimes_2 \rho_{B_n})= S(\rho_{A_n}) +
S(\rho_{B_n})- S(\rho_{A_nB_n}) \ .
$$
To simplify notations, let us write $S(A_n):= S(\rho_{A_n}),S(A_n\cup
B_n):= S(\rho_{A_nB_n}).$ Then we have
$$F(A,B) = \lim_{n\rightarrow \infty} S(A_n)+ S(B_n)- S(A_n\cup
B_n) \ .
$$
It follows that
$$
 F(A\cup B, A\cup C)= \lim_{n\rightarrow \infty} \big(S(A_n\cup B_n) + S(A_n\cup
 C_n) - S(A_n) -S(A_n\cup B_n\cup C_n)\big).
 $$
 Note that $$S(A_n\cup B_n) + S(A_n\cup
 C_n) - S(A_n) -S(A_n\cup B_n\cup C_n) \geq 0 $$
 by strong subadditivity of von Neumann entropy, (1) follows and (2)
 also follows from the limit formula and the fact that $F(A,B)$ is finite by Theorem
 \ref{main1}. \par

 (3): Starting  with arbitrary real valued function $G$ defined on
 open connected
 intervals of $S^1,$  we can define $G(A)$ for any $A\in \PI$ as
 follows: define
 $G(A\cup B) = G(A)+ G(B)- F(A, B)$
 when $A$ and $B$ are disjoint. It is easy to see that such
 $G(A\cup B)$ is well defined and only depends on $A\cup B$ thanks
 to (2). \par
 (4): This follows from Theorem \ref{main1}, (1) and direct
 computations. \par
 (5): The identities follow from (3).

 (6): We note that by  Theorem \ref{main1} and monotonicity of
 relative entropy in (4) of Th. \ref{515} that for $\B$,  $F_\B (A,B)\leq F_{\A_r} (A,B)<\infty.$ For (1) and (2) we
 can use remark \ref{invf} and proceed in exactly the same way as in free fermion net case.
 (3) and (5)  are proved in the same way as in free fermion net case.

 \qed

\subsection{Structure of singularities in the finite index
case}\label{structure}

$G$ from (3) in Th. \ref{soft1}  can be thought as ``regularized"
version of von Neumann entropy which is always infinite in our case
(cf. \cite{Nar}) . From (3) of the above Theorem we see that if we
only allow $G$ to be defined on $\PI$ then $G$ is highly non unique.
Due to the continuity properties of $F(A,B)$, we require that $G(A)$
 depends continuously only on the length $r_A$ of interval $A$. In
addition we require that $G(A)= G(A^c)$ for a connected interval,
and we set $G(\emptyset )=0$.  Still such $G$ is highly non unique.
However, we shall impose further conditions coming from studying the
singularities of relative entropy when we allow intervals to
approach each other. Let $B_\epsilon =(a_1, a_{2\epsilon}) \cup C
=(a_2, b_2)\in \PI,$ with $|a_{2\epsilon}-a_2|=\epsilon>0.$ We shall
consider the singular limit when $\epsilon$ goes to zero while
fixing $a_1$ and $C.$ Let $B_0=(a_1,a_2).$ We will denote by
$B_0\bar \cup C= (a_1,b_2),$ i.e., $B_0\bar \cup C$ is obtained from
$B_0 \cup C$ by adding the point $a_2$: notice in the process  the
number of components decrease by $1$.

To probe the singularity structure of  von Neumann entropy, we can
consider $F( B_\epsilon, C)$ which goes to $\infty$ as
$\epsilon\rightarrow 0$ while fixing $a_1$ and $C.$ As an example,
by Th. \ref{main1}
$$
F_{\A_r }(B_\epsilon, C) = \frac{r}{6} \big(\ln |a_2-a_1| + \ln
|b_2-a_2|- \ln |b_2-a_1| -\ln (\epsilon)\big) +o(\epsilon) \ .
$$
Since $G(B_\epsilon\cup C)= G(B_\epsilon)+G(C)-F( B_\epsilon, C),$
the singularity structure of $G(B_\epsilon\cup C)$ is the same as
that of $-F( B_\epsilon, C)$ as $\epsilon\rightarrow 0$. In fact
this is also true for general case: consider
$$
G(A\cup B_\epsilon \cup C) = G(A)+ G(B_\epsilon \cup C) - F(A,
B_\epsilon \cup C) \ .
$$
One can see that the singularity  structure of $G(A\cup
B_\epsilon\cup  C)$ is the same as that of $G(B_\epsilon \cup C)$ as
$\epsilon\rightarrow 0$, since the rest of terms are bounded. So we
can not expect $G( B_\epsilon \cup C)$ to be close to
$$
G(B_0\bar \cup C)
$$
when $\epsilon\rightarrow 0$, but we may demand that \ben\label{Gs}
\lim_{\epsilon\rightarrow 0}  G(B_\epsilon \cup C) -P(\epsilon)=
G(B_0\bar \cup C) \een
for some function $P(\epsilon)$ which is independent of $B,C$. The
equation is a condition that connects the value of $G$ for different
components and as we shall see is a very useful condition. Equation
(\ref{Gs}) is of course equivalent to \ben\label{Gs2}
 G(B_0\bar \cup
C)= G(B_0) + G(C) - \lim_{\epsilon\rightarrow 0}\big(P(\epsilon)+
F(B_\epsilon, C)\big) \ . \een
In general we may take multiple singular limits. Equation (\ref{Gs})
allows us to evaluate such limits. Let us consider such an example
in details. Let $A=(a_2,b_2)$, $B_{\epsilon_1}=
(a_1,a_{2\epsilon_1})$, $C_{\epsilon_2}= (b_{2\epsilon_2},
b_3),|a_{2\epsilon_1}- a_2|= \epsilon_1>0$,  $|b_{2\epsilon_2}- b_2|=
\epsilon_2>0 $. Let $\epsilon_1$ goes to $0$ first, we find
$$
F(A\bar \cup B_0, A\cup C_{\epsilon_2}) = G(A\bar \cup B_0)+G(A\cup
C_{\epsilon_2})- G(A\bar \cup B_0\cup C_{\epsilon_2}) -G(A)
$$
since the same function $P(\epsilon_1)$  appears in both
$G(A\cup B_{\epsilon_1})$ and $G(A \cup B_{\epsilon_1}\cup
C_{\epsilon_2})$ with opposite signs. Then let $\epsilon_1$ goes to
$0$ we get by the same argument
$$
F(A\bar \cup B_0, A\bar \cup C_0) = G(A\bar \cup B_0)+G(A\bar \cup
C_0)- G(A\bar \cup B_0\bar \cup C_0) -G(A) \ .
$$
It is easy to see that the result is independent of the order of
taking limits, and this way we can extend the definition of $F(A,B)$
to any $F(A,B)$ with $A\in \PI, B\in \PI.$ Such $F(A,B)$ is used in
\cite{Casc}. In the case of free fermions, by Th. \ref{main1} we
have  that $P(\epsilon)= r/6 \ln \epsilon +o(\epsilon),$ and  we
have
$$
F(A\bar \cup B_0, A\bar \cup C_0) = -\frac{r}{6} \ln
\left|\frac{(b_2-a_2)(b_3-a_1)}{(b_3-a_2)(b_2-a_1)}\right| \ .
$$
 Now we will show that equation (\ref{Gs}) is also true for a large
class of examples. We assume that $\B\subset \A_{r}$  has finite
index.

Note that by Lemma \ref{index} $\mu_\B = \mu_{\A_r} [\A: \B]^2= [\A:
\B]^2=\lambda^{-2} .$

 Let $F_1(A, B):= F_{\A_r} (A,B)- F_\B (A,B)$ and
$G_1(A)= G_{\A_r}(A)-G_\B(A)$.  Then $F_1(A, B)$ verifies (2) and (3)
of Th. \ref{soft1}. Note that $F_1(A, B)$ is not non-negative in
general, being the difference of two non-negative numbers, but is
always bounded by finite index assumptions.

We examine possible solutions of equation (\ref{Gs2}) for $G_1.$ Let
$B_\epsilon, C$ be two connected intervals as in  equation
(\ref{Gs2}), and $E$ the unique conditional expectation from
$\A_r(B_\epsilon) \vee \A_r(C)$ to $\B(B_\epsilon) \vee \B(C)$ which
preserves the state $\omega_1\otimes_2 \omega_2$. Then
$S_\A(\omega,\omega_1\otimes_2 \omega_2)
=S_\B(\omega,\omega_1\otimes_2 \omega_2) + S(\omega, \omega\cdot E)$
by Th. \ref{515}.  Note that by Pimsner-Popa inequality $E(x)\geq
\lambda^{-2}x$ for positive $x,$ and so $F_1(B_\epsilon,C)
=S(\omega, \omega\cdot E)\leq \ln \lambda^{-2}.$ By  Th. \ref{main2}
$\lim_{\epsilon \rightarrow 0}F_1(B_\epsilon,C)= \ln\lambda^{-1},$
 and equation (\ref{Gs2}) is simply
$$
 G_1(B_0\bar \cup
C)= G_1(B_0) + G_1(C) - (P-\ln \lambda) \ ,
$$
where $P$ is a constant. Up to a constant in the definition of
$G_1(A)$ we can set $P=\ln \lambda$, and it follows that  $G_1(A)$
is a constant multiplied by the arc length of $A$.   But since we
also require $G_1(A) = G_1(A^c)$, $G_1(A)=0.$

In this case we get $G_\B (A) = G_\A (A)$  for any connected
interval $A$, and use $G_\B= G_{\A_r}- G_1$ the system of solutions
of equation (\ref{Gs2}) for $\B$.

We have proved the following :

\bthm\label{soft2} Assume that a subnet $\B\subset \A_{r}$  has
finite index, then:

(1): $G_\B ((a,b))= \frac{r}{6} \ln |b-a|$ and verifies equation
(\ref{Gs2}) and (3) of Th. \ref{soft1}, and
$$
F_\B (A, B) = -\frac{r}{6} |\ln \eta_{AB}| \ ,
$$ where $A, B$ are two overlapping intervals with cross ratio
$0<\eta_{AB}<1$; \par

(2) Let $B=(a_1, a_{2\epsilon})$, $C=(a_2 ,b_2)$, $|a_{2\epsilon}-a_2|=
\epsilon >0$. Then:
$$
F_\B (B, C) = \frac{r}{6} \big(\ln |a_2-a_1| +\ln |b_2-a_2|- \ln
|b_2-a_1| -\ln (\epsilon)\big) -\frac{1}{2}\ln \mu_\B + o(\epsilon)
$$
as $\epsilon $ goes to $0$. \ethm

In exactly the same way if $\B\subset \C$ is a subnet with finite
index where $\B$ is as in the above theorem, then we also get a
system of solutions of equation (\ref{Gs2}) for $\C$ as in the above
theorem.

\brem\label{kitaevr}
 It is  interesting to note that  the
constant term in (2) of Th. \ref{soft2} seems to be related to the
topological entropy discussed in \cite{Kitaev} even with the right
factor: in our case we have additional factor $1/2$ since we are
discussing chiral half of CFT. \erem

We conjecture that the above theorem is true for any rational
conformal net, where $r$ is replaced by the central charge. More
examples where Th. \ref{soft2} applies are discussed in Section
\ref{examples}.

Notice also that the cross ratio enters in formulas concerning nuclearity (partition function)
\cite{BAL} and entanglement entropy \cite{Hollands}, so we can infer relations about
the mutual information and these quantities.

\subsubsection{Failure of duality is related to global
dimension}\label{fail}

By Th. \ref{main1} for the free fermion net $\A_r,$ and two intervals
$A=(a_1, b_1)$, $B=(a_2, b_2)$, where $b_1< a_2,$ we have
$$F_\A (A,B)=\frac{-r}{6}\ln \eta \ ,$$
where $\eta = \frac{(b_1-a_2)(b_2-a_1)}{(b_1-a_1)(b_2-a_2)}$ is the
cross ratio, $0< \eta <1.$ For simplicity we denote by $F_{\A_r}
(\eta)= F_\A (A,B)$.

One checks that  $F_{\A_r} (A,B)= F_{\A_r} (A^c,B^c)$, which is in
fact equivalent to
$$
F_{\A_r} (\eta)- F_{\A_r}(1-\eta)=\frac{-r}{6}\ln
\left(\frac{\eta}{1-\eta}\right) \ .
$$
Similarly for $\B\subset {\A_r}$ with finite index, by Th.
\ref{soft2} $F_\B(A,B)= F_\B (A^c,B^c)$ is equivalent to
$$
F_\B (\eta)- F_\B (1-\eta)=\frac{-r}{6}\ln \left(\frac{\eta}{1-\eta}\right) \ .
$$
 We note that $F_{\A_r} (A,B)= F_{\A_r} (A^c,B^c)$ for the free fermion net
$\A_r.$ However here we show that  $F_\B (A,B)\neq  F_\B (A^c,B^c)$
with $\B\subset {\A_r}$  has finite index $[{\A_r}:\B]=\lambda^{-1}
>1.$ By Lemma \ref{index} $\mu_\B=[{\A_r}:\B]^2. $

We note that, as before the proof of Th. \ref{soft2}, $S(\omega,
\omega\cdot E)= F_1 (\eta)=F_\A (\eta)- F_\B (\eta) $ is a
decreasing function of $\eta$, and $0\leq F_1 (\eta)\leq F_\A
(\eta).$ So we have
$$
\lim_{\eta\rightarrow 1}F_1 (\eta)=0 \ .
$$
On the other hand, by Th. \ref{main2}
$$
\lim_{\eta\rightarrow 0}F_1 (\eta)= \ln [\A_r:\B] = \frac{1}{2}
\ln\mu_\B \ .
$$
It follows that $F_\B (A,B)\neq  F_\B (A^c,B^c)$ due to the fact
that $\mu_\B>1$.

\subsection{Computation of limit of relative entropy}
In this section we determine the exact limit of relative entropies
which are necessary for analyzing the singularity structures of
entropies  in Section \ref{structure}. The goal is to prove the
following:

\bthm\label{main2} Assume that subnet $\B\subset \A$ has finite
index, $\B$ is strongly additive.  Let $I_1$ and $I_2$ be two
intervals obtained from an interval $I$ by removing an interior
point, and let $J_{n}\subset I_2, n\geq 1$ be an increasing sequence
of intervals such that
$$
\bigcup_n J_{n} =I_2,\quad \bar J_n\cap \bar I_1=\emptyset \ .
$$
Let $E_n$ be the conditional expectation from $ \A(I_1)\vee \A(J_n)$
to  $\A(I_1)\vee \B(J_n)$ such that $E_n(xy)= x E_I(y), \forall x\in
\A(I_1),y\in \A(J_n).$ Then
$$
\lim_{n\rightarrow \infty} S(\omega, \omega\cdot E_n)= [\A:\B]\ .
$$
\ethm

\subsubsection{Basic idea from Kosaki's formula}

Denote by $\phi_n= \omega\cdot E_n$.  By Kosaki's formula (cf.
\cite{Kos})
$$
 S(\omega, \omega\cdot E_n) =\sup_{m \in \mathbb{N}}\sup_{x_t+y_t=1}
 \left(\ln k - \int_{k^{-1}}^\infty\Big(\omega (x_t^*x_t)\frac{1}{t} + \phi_n
 (y_ty_t^*)\frac{1}{t^2}\Big)dt\right) \ ,
 $$
 where $x_t$ is a step function which is equal to $0$ when $t$ is
 sufficiently large. To motivate
 the proof of  Th. \ref{main2}, it is instructive to see how we can
 get  $S(\omega, \lambda\omega)= -\ln \lambda, 0<\lambda<1$ from
 Kosaki's formula. By tracing the proof in \cite{Kos}, one can see
 that the path which gives approximation to  $-\ln \lambda$ is given
 by the following continuous path
 $$
 x(t)= \frac{\lambda}{\lambda+t}, y(t)= \frac{t}{\lambda+t}, t\geq
 k^{-1}
 $$
 and with such a choice we have
$$\ln k - \int_{k^{-1}}^\infty\Big(\omega (x_t^*x_t)\frac{1}{t} + \phi_n
 (y_ty_t^*)\frac{1}{t^2}\Big)dt
 = -\ln  (\lambda+1/k)
 $$ which tends to  $-\ln \lambda$ as $k$ goes to $\infty$. This
 suggests that for the proof of Th. \ref{main2}, we need to choose
 path $x_t, y_t$ such that
$\omega (x_t^*x_t)$ and $ \phi_n
 (y_ty_t^*)$ are close to $\big(\frac{\lambda}{\lambda+t}\big)^2$ and
 $\lambda \big(\frac{t}{\lambda+t}\big)^2$ respectively, and this motivates
 our Prop. \ref{Kos1} and the proof of Th. \ref{main2}.

\subsubsection{A key step in the proof of  Th. \ref{main2}}

Let $e_1\in \A(I_1), e_2\in \A(J_1)$ be Jones projections for
$\B(I_1)\subset \A(I_1)$ and $\B(J_1)\subset \A(J_1)$ respectively
as in  Lemma \ref{jones}.  Let $P$  be the projection from the
vacuum representation of $\A$ onto the vacuum representation of
$\B$. By Lemma \ref{jones},  there is a unitary $u\in \B(I)$ such
that $ue_1u^*= e_2.$  Choose isometry $v_2\in \B(J_1)$ such that
$\lambda^{-1} v_2^* e_2 v_2 = 1.$ Note that $e_2 v_2v_2^* e_2=
\lambda e_2,$  and $Pe_2P =\lambda P.$ It follows that $Pe_2^+P
=\lambda P, Pe_2^-P =0$ by our assumption that  $[\Gamma, P]=0.$

Since $\B$ is strongly additive, we can find a sequence of bounded
operators $u_n \in \B(I_1)\vee \B(J_n), n\geq 2 $ such that
$u_n\rightarrow u$ strongly. Let $e_{2n}:= u_n e_1 u_n^*.$ Then
$e_{2n}\rightarrow e_2$ strongly.

\bprop\label{Kos1} For any $\epsilon>0$, one can find $n\geq 2$ and
$e\in \A(I_1)\vee \A(J_n)$ such that
$$
|\omega (e)-1|< \epsilon, |\omega (e^*)-1|< \epsilon,\ |\omega
(e^*e)-1|< \epsilon,\ |\phi_n (ee^*)-\lambda|< \epsilon \ .
$$\eprop

\proof

Let us first denote by $e= \lambda^{-1} v_2^* e_{2n} e_2 v_2^*.$ We
will show that given $\epsilon>0$, we can choose $n$ sufficiently
large such that $e$ verifies the conditions in the Proposition.
First we observe that since $e_{2n}\rightarrow e_2$ strongly, it
follows that $e\rightarrow 1$ strongly, and hence by choosing $n$
sufficiently large we can have
$$
|\omega (e)-1|< \epsilon, |\omega (e^*)-1|< \epsilon, |\omega
(e^*e)-1|< \epsilon
$$
Now let us evaluate
$$
\phi_n (ee^*)= \phi_n (\lambda^{-2} v_2^* e_{2n} e_2v_2^* v_2 e_2
e_{2n} v_2) = \lambda^{-1}\phi_n (v_2^* e_{2n}  e_2e_{2n} v_2)
$$
Recall the definition of $\phi_n$ as a state on $\A(I_1)\vee
\A(J_n)$: For any  $x,y$ with $x\in \A(I_1), y\in \A(J_n)$,
$$
\phi_n(xy) = \langle \Omega, xPyP\Omega\rangle
$$
Recall that $e_2= e_2^+ + e_2^-,Pe_2^+P =\lambda P, Pe_2^-P =0.$ To
evaluate $\phi_n (v_2^* u_n e_1 u_n^*   e_2   u_n e_1 u_n^*v_2),$ we
approximate  $u_n$ with finite linear combination of operator of the
form $u_{1m} u_{2m}$ with $u_{1m}\in \B(I_1),u_{2m}\in \B(J_n),$,
then we move those operators in $\A(I_1)$ to the left of those
operators in $\A(J_n)$ using commuting or anti-commuting relations,
and  it is crucial to observe  the operators that belong to
$\A(J_n)$ has only one term $e_2^+$ or  $e_2^-$£¬ and the rest are
in $\B(J_n).$ When compressed such term with $P$ and acting on
$\Omega,$ we see that $e_2^+$ is replaced with $\lambda,$ and
$e_2^-$ is replaced with $0$. We note that $e_2^+$ commuting with
$\A(I_1).$ It follows that
$$
\phi_n (ee^*)= \phi_n (\lambda^{-2} v_2^* e_{2n} e_2v_2^* v_2 e_2
e_{2n} v_2) = \lambda^{-1}\phi_n (v_2^* e_{2n}  e_2e_{2n} v_2) =
\langle \Omega, v_2^* (u_n e_1 u_n^*)^2 v_2 \Omega\rangle \ .
$$
Since $ v_2^* (u_n e_1 u_n^*)^2 v_2$ goes to $ v_2^* e_2
v_2=\lambda $ strongly, the Proposition is proved . \qed

\subsubsection{The proof of  Th. \ref{main2}} Recall $\phi_n
=\omega\cdot E_n.$ By Pimsner-Popa inequality, $E_n(x)\geq \lambda
x$ for any positive $ x\in \A(I_1)\vee \A(J_n),$ it follows that
$\phi_n \geq \lambda \omega$, and hence by Th. \ref{515}
$$ S(\omega, \omega\cdot E_n) \leq  [\A:\B].
$$
Note that by monotonicity of relative entropy $S(\omega, \omega\cdot
E_n)$ increases with $n,$ hence $\lim_{n\rightarrow \infty} S(\omega,
\omega\cdot E_n)$ exists and is less or equal to $[\A:\B].$

By Kosaki's formula
$$
 S(\omega, \omega\cdot E_n) =\sup_{m \in \mathbb{N}}\sup_{x_t+y_t=1}
 \left(\ln k - \int_{k^{-1}}^\infty\Big(\omega (x_t^*x_t)\frac{1}{t} + \phi_n
 (y_ty_t^*)\frac{1}{t^2}\Big)dt \right) \ ,
 $$
 where $x_t$ is a step function which is equal to $0$ when $t$ is
 sufficiently large. Since we can approximate any continuous
 function with step functions in the strong topology and vice versa,  we can assume that
 $x_t$ is continuous and   is equal to $0$ when $t$ is
 sufficiently large. Given $\epsilon >0$, for fixed $k,m \in
 \mathbb{N}$ choose $e$ as in Proposition \ref{Kos1} and $$
 x_t= 1- \frac{t}{\lambda +t} e, k^{-1} \leq t\leq m \ .
 $$
We have
$$
\omega(x_t^*x_t)= 1-\frac {t}{\lambda+t}\omega (e)-\frac
{t}{\lambda+t}\omega (e^*) +\left(\frac {t}{\lambda+t}\right)^2 \omega (e^*e)
$$
and
$$
\phi_n (y_ty_t^*)=\left(\frac {t}{\lambda+t}\right)^2\phi_n (ee^*) \ .
$$
By Proposition \ref{Kos1} we can choose $n$ large
 enough such that
 $$
 \int_{k^{-1}}^m \Big|\omega(x_tx_t^*) - \Big(\frac {\lambda}{\lambda+t}\Big)^2 \Big| \frac{dt}{t}\leq \epsilon \ ,
 $$
$$
 \int_{k^{-1}}^m \Big|\phi_n (y_ty_t^*) - \lambda\Big(\frac {t}{\lambda+t}\Big)^2\Big| \frac{dt}{t^2}\leq \epsilon \ ,
$$
and with such a choice of $n$ we have:
\begin{multline*}
\ln k - \int_{k^{-1}}^\infty \Big(\omega (x_t^*x_t)\frac{1}{t} + \phi_n
 (y_ty_t^*)\frac{1}{t^2}\Big)dt \geq  \\
\ln k - \int_{k^{-1}}^m \left(\Big(\frac {\lambda}{\lambda+t}\Big)^2 \frac{1}{t}
+\Big(\frac {t}{\lambda+t}\Big)^2\frac{\lambda}{t^2}\right)dt +1/m  -
2\epsilon \\
= \ln \Big(\frac{k}{k\lambda +1}\Big) -\ln \Big(\frac{m}{\lambda+m}\Big) +1/m
-2\epsilon \ .
\end{multline*}
It follows that
$$
\lim_{n\rightarrow \infty} S(\omega, \omega\cdot E_n) \geq \ln
\Big(\frac{k}{k\lambda +1}\Big) -\ln \Big(\frac{m}{\lambda+m}\Big) +1/m -2\epsilon \ .
$$
Let $k,m$ go to $\infty$ and $\epsilon$ go to $0$, we have proved
theorem. \qed

\subsection{More Examples}\label{examples}
\subsubsection{Orbifold examples}
Take $U(1)_{4k^2}\subset U(1)_1.$  This is $\mathbb{Z}_{2k}$
orbifold of $U(1)_1.$ So Th. \ref{soft2} apply to the net
$U(1)_{4k^2}$. Another special case is when $k=1,$ we can take a
further $\mathbb{Z}_2$ orbifold of  $U(1)_{4}$ which corresponds to
complex conjugation on $U(1)$  to get a tensor product of two Ising
model with central charge $\frac{1}{2}.$ It follows that Ising model
with central charge $\frac{1}{2}$ verifies Th. \ref{soft2}, and in
particular violates duality discussed in Section \ref{fail}.\par

More generally, we can take any finite subgroup of $U(n)$ which
commutes with $\Ad \Gamma$ and obtain orbifold subnet of $U(n)_1.$
This provides  a large family of examples which verify Th.
\ref{soft2}.
\subsubsection{Conformal inclusions}
By \cite{XN}, we have the following inclusions with finite index:
$$
SU(n)_m \times SU(m)_n \times U(1)_{mn(m+n)^2}\subset
Spin(2mn)_1\subset U(mn)_1 \ .
$$
So  Th. \ref{soft2} apply to the net $SU(n)_m \times SU(m)_n \times
U(1)_{mn(m+n)^2}.$ If we take $m=n$, then since $U(1)_{(4n^4)}$
verifies  Th. \ref{soft2} by the example in  previous section, it
follows that the net associated with $SU(n)_n \times SU(n)_n$, and
hence  the net associated with $SU(n)_n$ also verifies  Th.
\ref{soft2}.

\vspace{0.4cm}

\noindent
{\bf Acknowledgements}. The authors would like to thank E. Witten for stimulating, enlightening
discussions and encouragement. We also thank Y.  Tanimoto for comments.

{\footnotesize

\end{document}